\documentclass{IOS-Book-Article}

\usepackage{times}
\normalfont
\usepackage[T1]{fontenc}

\usepackage{amssymb}
\usepackage{amsmath}
\usepackage{algorithmic}
\usepackage{graphicx}
\newcommand{\Gr}{Gr\"obner }
\newcommand{\cL}{\cal{L}}

\newcommand{\Q}{\mathbb{Q}}

\newcommand{\K}{\mathbb{K}}
\newcommand{\N}{\mathbb{N}}
\newcommand{\M}{\mathbb{M}}
\newcommand{\R}{\mathbb{R}}
\newcommand{\X}{\mathbb{X}}

\newcommand{\lcm}{\mathop{\mathrm{lcm}}\nolimits}
\newcommand{\lm}{\mathop{\mathrm{lm}}\nolimits}
\newcommand{\lt}{\mathop{\mathrm{lt}}\nolimits}

\newcommand{\pol}{\mathop{\mathrm{pol}}\nolimits}
\newcommand{\anc}{\mathop{\mathrm{anc}}\nolimits}
\newcommand{\nmp}{\mathop{\mathrm{nmp}}\nolimits}

\newcommand{\Id}{\mathop{\mathrm{Id}}\nolimits}
\newcommand{\idx}{\mathop{\mathrm{idx}}\nolimits}

\def \bg #1 {\begin{tabular}{{#1}}}
\def \nd {\end{tabular}}

\renewenvironment{algorithm}[1]{
  \begin{center}
    {\bf Algorithm: #1}\\*
     \begin{tabular}{|p{120mm}|} \hline
} {
 \\ \hline
 \end{tabular}
 \end{center}
}

\begin{document}
\begin{frontmatter}

\title{Involutive Algorithms for Computing \Gr Bases}
\runningtitle{Involutive Algorithms}

\author{\fnms{Vladimir P.} \snm{Gerdt}
\thanks{Vladimir P. Gerdt, 141980 Dubna, Russia.
E-mail: gerdt@jinr.ru.}} \runningauthor{V.P. Gerdt}

\address{Laboratory of Information Technologies\\
        Joint Institute for Nuclear Research}


\begin{abstract}
In this paper we describe an efficient involutive algorithm
for constructing \Gr bases of polynomial ideals. The algorithm
is based on the concept of involutive monomial division which restricts
the conventional division in a certain way. In the presented algorithm
a reduced \Gr basis is the internally fixed subset of an involutive
basis, and having computed the later, the former can be output
without any extra computational costs. We also discuss some accounts of
experimental superiority of the involutive algorithm over Buchberger's
algorithm.
\end{abstract}
\end{frontmatter}


\section{Introduction}

The concept of \Gr bases for polynomial ideals together with the algorithmic characterization of
these bases in terms of $S-$polynomials was invented almost 40 years ago by Buchberger in his
PhD thesis~\cite{BB65}. Since that time \Gr bases have become the most universal
algorithmic tool in commutative algebra and algebraic geometry~(see, for example,
books \cite{BW}-\cite{GP}). Moreover, this role of \Gr bases was
extended to rings of linear differential and shift operators and to some other
noncommutative algebras~(see~\cite{BW98,CAHandbook} and references therein).

The fundamental property of a \Gr basis of a polynomial ideal is the divisibility of any element
in the related initial ideal by the leading monomial of an element in the basis. In our
paper~\cite{GB98a} an algorithmic approach was developed based on an axiomatically
formulated restriction for the conventional monomial division. This restricted division which
we called involutive generalizes the pioneering results by Zharkov and Blinkov~\cite{ZB} on
carrying over the involutive methods from differential equations to commutative algebra and
leads to a more general concept of involutive bases. An involutive basis is a (generally redundant)
\Gr basis. But for all that, an element in the initial ideal is multiple of the leading term of an
element in the involutive basis not only in the conventional sense but also in the involutive
sense. Just as a monic reduced \Gr basis~\cite{Buch85}, a monic minimal involutive basis~\cite{GB98b}
of an ideal is uniquely defined for a given monomial order.

An impetus was given to the involutive division technique with the experimental
demonstration in~\cite{GB98a,ZB} the efficiency of involutive algorithms as an alternative
method to compute \Gr bases. The concepts, ideas and algorithms developed in~\cite{GB98a} got
further development and extension to linear differential algebra and noncommutative rings in our
papers~\cite{GB98b}-\cite{G02} and in the papers of other authors~\cite{CG01}-\cite{Evans}.

\noindent
With all this going on, the involutive algorithms have been implemented for particular
involutive divisions in Reduce~\cite{GB98a,ZB,GBY01}, Mathematica~\cite{Berth},
Maple~\cite{Daniel}, MuPAD~\cite{SH}, C/C++~\cite{GBY01} and recently in Singular
(version 2.0.5).

Having analyzed our paper~\cite{GB98a}, Apel~\cite{Apel} introduced a slightly different
concept of involutive division and designed another form of involutive algorithms based on
his concept. The approach of Apel was further elaborated by himself and
Hemmecke~\cite{Hem02}-\cite{Hem03}. Differences of both approaches~\cite{GB98a}
and~\cite{Apel} and some unifying aspects were studied in~\cite{HemThesis}.

For a given finite polynomial set, an involutive division partitions the variables
into two disjoint subsets called multiplicative and non-multiplicative. The construction
of an involutive basis is often called completion to involution and consists in
examining the involutive reducibility of non-multiplicative prolongations, i.e.,
products of the polynomials and their non-multiplicative variables. In doing so,
the involutive algorithms also examine all essential $S-$polynomials~\cite{GB98b,Apel}
but with rather special reductions based on involutive division. Thereby, there is no
conceptual distinction between Buchberger algorithm and involutive algorithms. The
difference is in both critical pairs selection and their reduction provided with the
presence of extra polynomials at involutive completion which are redundant in the
conventional sense.

Computer experiments with our C/C++ code~\cite{GBY01} implementing the involutive algorithm
for Janet division and with the standard benchmarks~\cite{BM,Verschelde} show that in the
most cases our code is faster than the most optimized modern-day implementations of
Buchberger's algorithm~\footnote{See Web site {\tt http://invo.jinr.ru} for comparison
of our involutive code with implementations of Buchberger's algorithm in Singular and Magma.}.
Thus the involutive algorithm may be considered as an improvement of Buchberger's algorithm.

In the present paper we follow the concept of involutive division in paper~\cite{GB98a},
which is briefly described in Sect.3 as a particular case of the restricted monomial division
introduced in Sect.2, and present our involutive algorithm in its detailed and optimized form
(Sect.4). The algorithm can also be used to output the reduced
\Gr bases of polynomial ideals without extra reduction of the computed involutive basis.
Some aspects of computational efficiency of the involutive algorithm are considered in Sect.5.
It is easy to rewrite the algorithm for linear differential polynomial ideals and for ideals
in rings of noncommutative polynomials with algebra of solvable type. Besides, the algorithm
admits a straightforward extension to modules in these rings.

\section{Preliminaries}

In this paper we shall use the following notations and conventions:
\begin{itemize}

\item $\N_{\geq 0}$ is the set of nonnegative integers.

\item $\X=\{x_1,\ldots,x_n\}$ is the set of polynomial variables.

\item $X\subseteq \X$ is a subset (possibly empty) of the variables.

\item $\R=\K[\X]$ is a polynomial ring over a zero characteristic field $\K$.

\item $\Id(F)$ is the ideal in $\R$ generated by $F\subset \R$.

\item $\M=\{x_1^{i_1}\cdots x_n^{i_n} \mid i_k\in \N_{\geq 0},\
1\leq k\leq n\}$ is the monoid of monomials in $\R$.

$\M_X\subseteq \M$ is a submonoid of $\M$ generated by the power products
of variables in $X$.

\item $\deg_i(u)$ is the degree of $x_i$ in $u\in \M$.

\item $\deg(u)=\sum_{i=1}^n \deg_i(u)$ is the total degree of $u$.

\item $\succ$ is an admissible monomial order such that
$x_1\succ x_2\succ\cdots\succ x_n$.

\item $u\mid v$ is the conventional divisibility relation of
monomial $v$ by monomial $u$.
If $u\mid v$ and 
$\deg(u)<\deg(v)$, i.e. $u$ is a proper
divisor of $v$, we shall write $u\sqsubset v$.

\item $\lm(f)$ and $\lt(f)$ are the leading monomial and the
leading term of $f\in
\R\setminus \{0\}$.

\item $\lm(F)$ is the leading monomial set for $F\subset \R\setminus \{0\}$.

\item ${\cL}$ is an involutive monomial division.

\item $u\mid_{\cL} v$ is the divisibility relation of monomial $v$
by monomial $u$ defined
by involutive 
division $\cL$.

\item ${\cal{C}}(U)$ is the cone generated by monomial set $U\subset \M$.

\item ${\cal{C}}_{\cL}(U)$ is the ${\cL}-$cone generated by monomial set $U\subset \M$.

\item $NM_{\cL}(u,U)$ is the set of ${\cL}-$non-multiplicative
variables of
monomial $u\in U\subset \M$.

\item $HNF_{\cL}(f,F)$ is the ${\cL}-$head normal form of $f\in
\R$ modulo
$F\subset \R$.

\item $NF_{\cL}(f,F)$ is the ${\cL}-$(full) normal form of $f\in
\R$ modulo
$F\subset \R$.

\item $\pol(p)=f$ is the first element in triple $p=\{f,g,X\}$
where
 $f,g\in F\subset \R$. 

\item $\anc(p=\{f,g,X\})=g$ is
 an ancestor of $f$ in $F$.

\item $\nmp(p=\{f,g,X\})=X$ is a subset
of non-multiplicative variables for $f$. 

\item $\idx(f,F)$ is the index of an element $f\in F$ in an indexed set $F$.

\end{itemize}

\vskip 0.3cm
\noindent
To introduce the concept of involutive division and describe the involutive algorithm
for construction of \Gr bases we need the following definitions.

\vskip 0.3cm
\noindent
{\bf Definition 2.1.} \ \ A linear monomial order $\succ$ is called {\em admissible}
if the conditions
$$m\neq 1\ \Longleftrightarrow\ m\succ 1,\quad \ m_1\succ m_2\
 \Longleftrightarrow m_1m\ \succ m_2m$$
hold for any monomials $m,m_1,m_2\in \M$ .

\vskip 0.3cm
\noindent
Hereafter we consider only admissible monomial orders and omit the word ``admissible''.
\vskip 0.3cm
\noindent
{\bf Definition 2.2.}\ \ For a finite set $F\subset \R\setminus \{0\}$, a polynomial
$h\in \R\setminus \{0\}$ and a monomial order $\succ$ consider polynomial $\tilde{h}$ given by
$$
\tilde{h}=h-\sum_{ij}\alpha_{ij}m_{ij}f_j\,,
$$
where $\alpha_{ij}\in \K,\ f_j\in F,\ m_{ij}\in \M$ and $\lm(m_{ij}f_j)\preceq \lm(h)$ whenever
$\ \alpha_{ij}\neq 0$. If $\lm(\tilde{h})$ has no divisors in $\lm{F}$, then $\tilde{h}$
is called {\em a head normal form of\, $h$ modulo} $F$. In this case we shall write
$\tilde{h}=HNF(h,F)$. In addition, if all other monomials occurring in $\tilde{h}$
also have no divisors in $\lm(F)$, then $\tilde{h}$ is called {\em a normal form of\, $h$
modulo} $F$. In such a case we shall write $\tilde{h}=NF(h,F)$.

\vskip 0.3cm
\noindent
{\bf Definition 2.3.} \ \ Given an ideal ${\cal{I}}\subset \R$ and an order $\succ$, a finite subset
$G\subset \R$ is called {\em \Gr basis} of ${\cal{I}}$ if
\begin{equation}
(\forall f\in {\cal{I}})\ (\exists g\in G)\ \ [\ \lm(g)\mid \lm(f)\ ]\,, \label{conv_div}
\end{equation}
where $u \mid v$ denotes the conventional divisibility of monomial $v$ by monomial $u$.

\vskip 0.3cm
\noindent
From Definition 2.2 it follows that (\ref{conv_div}) is equivalent to
\begin{equation}
(\forall f\in I)\ [\ NF(f,G)=0\ ]\,. \label{conv_imp}
\end{equation}

\vskip 0.3cm
\noindent
Definitions 2.1-2.3 are well-known
ingredients of the \Gr bases theory~\cite{BB65}-\cite{BW98}. Now
instead of the conventional monomial division in~(\ref{conv_div})
we consider another divisibility relation and the related concepts
of the normal form and \Gr basis.

\vskip 0.3cm
\noindent
{\bf Definition 2.4.} \ \ A {\em restricted division} $r$ on $\M$ is a
transitive relation $u \mid_r v$\ ($u,v\in \M$) such that
$$ u \mid_r v \Longrightarrow u \mid v\,. $$
If $u \mid_r v$, then $u$ is {\em $r-$divisor}\, of $v$ and $v$ is {\em $r$-multiple} of $u$,
respectively. Note, that the whole class of restricted divisions includes the conventional
division as well.

\vskip 0.3cm
\noindent
{\bf Definition 2.5.} \ \ Let $\succ$ be a monomial order, $r$ be a restricted division
and $F\subset \R\setminus\{0\}$ be a finite polynomial set. If $g\in \R\setminus\{0\}$
has a monomial which is $r-$multiple of an element in $\lm(F)$, then $g$ is called
{\em $r-$reducible modulo $F$}. If $\lm(g)$ has an $r-$divisor in $\lm(F)$, polynomial $g$
is called {\em $r-$head reducible modulo $F$}. If
for any $f\in F$ and for any monomial $u$ occurring in $f$ there is no $r$-divisors of
$u$ in $\lm(F)\setminus \{lm(f)\}$, then $F$ is called {\em $r$-autoreduced}.

\vskip 0.3cm
\noindent
{\bf Definition 2.6.}\ \ Similarly to the conventional division in Definition 2.2,
an {\em $r-$head normal form} $HNF_r(h,F)$ and {\em $r-$normal form} $NF_r(h,F)$ of
$h\in \R\setminus \{0\}$ modulo $F\subset \R\setminus \{0\}$,
for a given restricted division $r$, is the polynomial $\tilde{h}$ given by
$$
\tilde{h}=h-\sum_{ij}\alpha_{ij}m_{ij}f_j\,.
$$
Here
$\lm(f) \mid_r \lm(f)\,m_{ij}$ if $\alpha_{ij}\neq 0$,  $\lm(m_{ij}f_j)\preceq \lm(h)$.
And, respectively, $\lm{h}$ and all monomials in $\tilde{h}$
are $r-$irreducible modulo $\lm(F)$.

\vskip 0.3cm
\noindent
{\bf Definition 2.7.} \ \ Given an ideal ${\cal{I}}\subset {\R}$,  an order $\succ$
and a restricted division $r$, a finite subset $G\subset \R$ is called {\em $r-$basis}
of ${\cal{I}}$ if
$$ (\forall f\in {\cal{I}})\ (\exists g\in G)\ \
 [\ \lm(g)\mid_r \lm(f)\ ]. $$

\vskip 0.3cm
\noindent
As well as in the case (\ref{conv_imp}) of the conventional monomial division,
$r-$basis $G$ can also be defined by the relation
$$(\forall f\in I)\ [\ NF_r(f,G)=0\ ]\,. $$
From Definition 2.4 it follows that a $r-$basis is always a \Gr basis.

\section{Involutive Division and Involutive Bases -- Definition}

A natural way to introduce a restricted monomial division $r$ is to indicate a certain
subset $X(u)\subseteq \X$ of variables for a monomial $u\in \M$ and to define for $v\in \M$
$$u \mid_r v \Longleftrightarrow v=u\cdot w, \quad w\in \M_{X(u)}\,,$$
where $\M_{X(u)}$ is the monoid of power products constructed from the variables in $X(u)$.

\noindent
Definitions 2.5 and 2.6 deal with $r$-divisors
taken from a fixed finite monomial set. By this reason it suffices to define an $r$-division
for an arbitrary finite set of divisors. The below definitions, taken from~\cite{GB98a},
allowed to introduce a wide class of restricted divisions providing an algorithmic
way for construction of $r-$bases~\cite{GB98b,G02,CG01}.

\vskip 0.3cm
\noindent
{\bf Definition 3.1.} \ \ We say that an {\em involutive division} $\cL$
is defined on $\M$ if for any nonempty finite monomial set $U\subset \M$ and for any
$u\in U$ there is defined a subset $M_{\cL}(u,U)\subseteq \X$ of variables generating
monoid ${\cL}(u,U)\equiv \M_{M_{\cL}(u,U)}$ such that the following
conditions hold
\begin{enumerate}
\item $v\in U\ \wedge \ u{\cL}(u,U)\cap v{\cL}(v,U)\not=\emptyset \Longrightarrow
 u\in v{\cL}(v,U)\ \vee \ v\in u{\cL}(u,U)$.
\item  $v\in U\ \wedge \ v\in u{\cL}(u,U)\Longrightarrow {\cL}(v,U)\subseteq {\cL}(u,U)$.
\item  $u\in V  \wedge \ V\subseteq U\Longrightarrow {\cL}(u,U)\subseteq {\cL}(u,V)$.
\end{enumerate}
Variables in $M_{\cL}(u,U)$ are called $({\cL}-)$\,{\em multiplicative} for $u$ and those in
$NM_{\cL}(u,U)\equiv \X\setminus M_{\cL}(u,U)$ are called (${\cL-}$)\,{\em non-multiplicative} for $u$,
respectively. If $w\in u{\cL}(u,U)$, then $u$ is called {\em ${\cL}-$(involutive) divisor}
of $w$. In this case we shall often write $u \mid_{\cL} w$.

\vskip 0.3cm
\noindent
The restricted division given by Definition 3.1 called involutive since it generalizes
the main properties of three different separations of variables used by Janet~\cite{Janet},
Thomas~\cite{Thomas} and Pommaret~\cite{Pommaret} in the algebraic analysis of partial
differential equation systems based on their completion to involution. Those
separations satisfy~\cite{GB98a} properties 1-3 in Definition 3.1 and, hence, generate
involutive divisions.

Now, as particular examples of involutive divisions, we present Janet and Pommaret
divisions~\cite{GB98a} which are used most. In addition, we consider the example of
involutive division introduced in~\cite{G02}. In doing so,
we indicate multiplicative or non-multiplicative variables only. The remaining variables
are to be considered as non-multiplicative or multiplicative, respectively.

\vskip 0.3cm
\noindent
{\bf Example 3.1.}\ ( {\em Janet division} ).
For each $1\leq i\leq n$ divide $U\in \M$ into groups labeled by non-negative integers
$d_1,\ldots,d_i$
$$
[d_1,\ldots,d_i]=\{\ u\ \in U\ |\ d_j=\deg_j(u),\ 1\leq j\leq i\ \}.
$$
$x_1$ is multiplicative
for $u\in U$ if $\deg_1(u)=\max\{\deg_1(v)\ |\ v\in U\}$.
For $i>1$ $x_i$ is multiplicative
for $u\in [d_1,\ldots,d_{i-1}]$ when
$\deg_i(u)=\max\{\deg_i(v)\ |\ v\in
[d_1,\ldots,d_{i-1}]\}.$

\vskip 0.3cm
\noindent
{\bf Example 3.2.}\ ( {\em Pommaret division} ). For
$v=x_1^{d_1}\cdots x_{k-1}^{d_{k-1}}x_k^{d_k}$ with $d_k>0$
variables $x_j,j\geq k$ are multiplicative.
For $v=1$ all the variables
are multiplicative.

\vskip 0.3cm
\noindent
{\bf Example 3.3.}\  ( {\em Lexicographically induced division} ).  A variable $x_i$ is
non-multiplicative for $u\in U$ if there is $v\in U$ such that
 $v\prec_{Lex} u$ and $deg_i(u)<deg_i(v)$, where $\succ_{Lex}$ denotes
 the pure lexicographical monomial order.

\vskip 0.3cm
\noindent
In the last example one can replace the lexicographical order with any other
(admissible) order to obtain another involutive division~\cite{G02}. The separation
of variables into multiplicative ($M$) and non-multiplicative ($NM$) for divisions
in Examples 3.1-3.3 and for the set of three monomials in three variables shown in Table 1.


\begin{table}[h]
\label{Divisions} \caption{Separations of variables for
$U=\{x_1^2x_3,x_1x_2,x_1x_3^2\}$}
\begin{center}
\bg {|c|c|c|c|c|c|c|} \hline\hline
Element & \multicolumn{6}{c|}{Involutive division}
\\ \cline{2-7}
in $U$ & \multicolumn{2}{c|}{Janet} &\multicolumn{2}{c|}{Pommaret} &
\multicolumn{2}{c|}{Lex. induced} \\ \cline{2-7}
 & $M$ & $NM$  & $M$   & $NM$   & $M$  & $NM$  \\
\hline
$x_1^2x_3$ & $x_1,x_2,x_3$ & $-$ & $x_3$ & $x_1,x_2$ & $x_1$ & $x_2,x_3$ \\
$x_1x_2$ & $x_2,x_3$ & $x_1$ & $x_2,x_3$ & $x_1$ & $x_1,x_2$ & $x_3$ \\
$x_1x_3^2$ & $x_3$ & $x_1,x_2$ & $x_2,x_3$ & $x_1$ & $x_1,x_2,x_3$ & $-$ \\
\hline \hline
\nd
\end{center}
\end{table}
\noindent
As it was shown for the first time by Zharkov and Blinkov~\cite{ZB},
completion of multivariate polynomial systems to involution gives an algorithm for
construction of \Gr bases. Zharkov and Blinkov used
Pommaret separation of variables and thus constructed what is called now
Pommaret bases. They proved termination of their completion
algorithm for zero-dimensional ideals whereas for positive-dimensional ideals
Pommaret bases may not exist (as finite sets), and, thus, the algorithm may not
terminate.

Termination properties of completion algorithms are determined
by the underlying involutive divisions and can be studied by                                                                     completion
of monomial sets~\cite{GB98a}. The following two definitions
elucidate this important aspect of involutive algorithms.

\vskip 0.3cm
\noindent
{\bf Definition 3.2.}\ \ Set ${\cal{C}}(U)=\cup_{u\in U}\,u\,\M$ is called the {\em
cone} generated by finite monomial set $U\subset \M$ and set
${\cal{C}}_{\cL}(U)=\cup_{u\in U}\,u\,{\cal{L}}(u,U)$ is called
{\em ${\cL}-$(involutive) cone} of $U$.

\vskip 0.3cm
\noindent
{\bf Definition 3.3.}\ \ A finite monomial set $\bar{U}\supseteq U$ is called
{\em ${\cL}-$completion} of
set $U$ if ${\cal{C}}(U)={\cal{C}}_{\cL}(\bar{U})$.
If
\begin{equation}
{\cal{C}}_{\cL}(U)={\cal{C}}(U) \label{completeness}
\end{equation}
set $U$ is called {\em ${\cL}$-complete} or {\em
${\cL}$-involutive}.
If every finite set $U$ admits ${\cL}-$completion, then involutive division ${\cL}$
is called {\em Noetherian}.

\vskip 0.3cm
\noindent
Janet division and Lexicographically induced division are Noetherian~\cite{GB98a,G02}.
As to Pommaret division, it is non-Noetherian~\cite{GB98a} what can be also seen by explicit
completion of the monomial set in Table 1. The completion can be performed either by the
special monomial algorithms~\cite{GB98a,GBY01} or by the general polynomial algorithms
described below. The resulting sets are given by:
\begin{eqnarray*}
&& \mbox{Janet}:\ \{x_1^2x_3,x_1x_2,x_1x_3^2,x_1^2x_2\}, \\
&& \mbox{Pommaret}:\ \{x_1^2x_3,x_1x_2,x_1x_3^2,x_1^2x_2,\ldots,x_1^{i+3}x_2,
\ldots,x_1^{j+3}x_3,\ldots\}\ i,j\in \N_{\geq 0},\\
&& \mbox{Lex. induced}:\ \{x_1^2x_3,x_1x_2,x_1x_3^2,x_1x_2x_3\}.
\end{eqnarray*}
Apart from Noetherianity, algorithmically ``good'' involutive divisions must be
{\em continuous} and
{\em constructive} in accordance with two definitions that follow.

\vskip 0.3cm
\noindent
{\bf Definition 3.4.}\ \ A set $U$ is called {\em locally ${\cL}-$involutive}
or {\em locally
${\cL}-$complete}
if \begin{equation}
(\forall u\in U)\ (\forall x\in NM_{\cal{L}}(u,U))\ \ [\ u\cdot x\in {\cal{C}}_{\cL}(U)\ ]\,.
\label{cont}
\end{equation}
An involutive division ${\cL}$ is called {\em continuous} if conditions (\ref{cont}) imply
equality (\ref{completeness}) for any monomial set $U\in \M$.

\vskip 0.3cm
\noindent
Unlike the condition~(\ref{completeness}) of involutivity/completeness, the
conditions~(\ref{cont}) of local
involutivity/completeness admit an algorithmic verification and are fundamental for
involutive algorithms. In paper~\cite{GB98a} we found a criterion which allows to check
continuity of an involutive division.

Given a monomial order $\succ$ and an involutive division ${\cL}$, the separation of
variables for
elements in a finite polynomial set $F\subset \R\setminus \{0\}$ is defined in terms
of the leading
monomial set $\lm(F)$. Now we can introduce the concept of involutive basis as an $r-$basis
in
Definition 2.6 specified for an involutive division. But for all that we follow the
definition
in~\cite{GB98a} and demand for involutive bases are to be autoreduced in accordance
with Definition 2.7.

\vskip 0.3cm
\noindent
{\bf Definition 3.5.} \ \ Let ${\cal{I}}\subset \R$ be a nonzero ideal,
${\cL}$ be an involutive division and $\succ$ be a monomial order. Then a finite
${\cL}$-autoreduced subset $G\subset \R$ is called (${\cL}-$){\em involutive basis}
of ${\cal{I}}$ if
\begin{equation}
(\,\forall f\in {\cal{I}}\,)\ (\,\exists \, g\in G\,) \ [\ \lm(g)\mid_{\cL} \lm(f)\ ]\,. \label{inv_bas}
\end{equation}
If division ${\cL}$ is continuous, then conditions (\ref{inv_bas}) are equivalent~\cite{GB98a} to the
following {\em involutivity conditions}
\begin{equation}
(\ \forall f\in G\ )\ (\ \forall x_i\in NM_{\cL}(\lm(f),\lm(G))\ )\
 \ [\ NF_{\cL}(x_i\cdot f,G)=0\ ]\,. \label{inv_bas_1}
\end{equation}
Here $NF_{\cL}(f,F)$ denotes the ${\cL}$-normal form of $f$ modulo $F$ according to Definition 2.6
where $r-$division replaced by ${\cL}$-division.

\vskip 0.3cm
\noindent
{\bf Definition 3.6.}\ \ The product of a polynomial by its non-multiplicative (multiplicative)
variable is called {\em non-multiplicative (multiplicative) prolongation} of the polynomial.

\vskip 0.3cm
\noindent
Involutivity conditions (\ref{inv_bas_1}) give an algorithmic characterization of involutive bases
much like the algorithmic characterization of \Gr bases by $S-$polynomials established in~\cite{BB65}.
One can say that involutive bases are characterized by ${\cL}-$reducibility to zero of all the
non-multiplicative prolongations of elements in the basis.

The basic idea behind an algorithmic construction of involutive bases is to check the involutivity
conditions. If they are not satisfied, then the nonzero ${\cL}-$normal forms are added
in a certain order to the system until all the involutivity conditions (\ref{inv_bas_1}) satisfied.
With all this going on, an algorithmically ``good'' involutive division, in addition to its
Noetherianity  and continuity, should satisfy the {\em constructivity properties}~\cite{GB98a}
in accordance to the definition:
\vskip 0.3cm
\noindent
{\bf Definition 3.7.}\ \ An involutive division ${\cL}$ is called {\em constructive}
if for any $U\subset \M$, $u\in U$, $x_i\in NM_{\cal{L}}(u,U)$ such that $u\cdot x_i\not \in
{\cal{C}}_{\cL}(U)$
and
$$
 (\forall v\in U)\ (\forall x_j\in NM_{\cL}(v,U))\ (\ v\cdot x_j \sqsubset u\cdot x_i\ )\ \
 [\ v\cdot x_j\in {\cal{C}}_{\cL}(U)\ ]
$$
the condition
$ (\forall w\in {\cal{C}}_{\cL}(U))
 \ \ [\ u\cdot x_i\not \in w{\cL}(w,U\cup \{w\})\ ]$ holds.
\vskip 0.3cm
\noindent
Constructivity ensures that in the course of the algorithm there no needs to enlarge
the intermediate basis with multiplicative prolongations and only non-multiplicative prolongations
must be examined for the enlargement (completion). Note that all three divisions of Examples 3.1-3.3
are continuous and constructive~\cite{GB98a,G02}.

As any $r-$basis, an involutive basis is a \Gr basis. ${\cL}$-reducibility implies
the conventional reducibility (i.e. reducibility with respect to the conventional division).
But the converse is not true in general. By this reason, an involutive basis is generally
redundant as the \Gr one. Moreover, a monic reduced \Gr basis is uniquely defined by an
ideal and a monomial order~\cite{BB65,Buch85} whereas this is not true for involutive bases as
shows the following simple bivariate example~\cite{GB98a}.

\vskip 0.3cm
\noindent
{\bf Example 3.4.}\ \ Consider ideal in $\Q[x,y]$ generated by $F=\{x^2y-1,xy^2-1\}$. For the
lexicographical order with $x\succ y$ the polynomial sets
\begin{eqnarray*}
G_1&=&\{xy^2-1,xy-y^2,x-y,y^3-1\}\,, \\[0.1cm]
G_2&=&\{x^2y-1,x^2-y^2,xy^3-y,xy^2-1,xy-y^2,x-y,y^3-1\}\,
\end{eqnarray*}
are both Janet bases of $\Id(F)$.
\vskip 0.3cm
\noindent
However, for a constructive division one can define~\cite{GB98b} a minimal involutive basis
which, similarly to a reduced \Gr basis, being monic is uniquely defined by an ideal and a monomial
order. If $G$ is such a monic minimal involutive basis, then for any other monic involutive basis
$G_1$ the inclusion $G\subset G_1$ holds. For Example 3.4. the minimal Janet basis is $\{x-y,y^3-1\}$
and coincides with the reduced \Gr basis.

\section{Involutive Bases---Construction}
In the rest of the paper we assume that the input involutive division in the below algorithms
is Noetherian, continuous and constructive. Our goal is to construct
minimal bases. Having this in mind, we shall often omit the word ``minimal''.

First, we present the simplest version of an algorithm for constructing involutive polynomial
bases and illustrate its work by the bivariate polynomial set
from Example 3.4. Then we describe an improved version of the algorithm.

In the below algorithm {\bf InvolutiveBasis I} the whole polynomial data are
partitioned into two subsets $G$ and $Q$. Set $G$ contains a part of the
intermediate basis. Another part of the intermediate basis is contained in set
$Q$ and includes also all the non-multiplicative prolongations of
polynomials in $G$ which must be examined in accordance to the
involutivity conditions (\ref{inv_bas_1}).

At the initialization step of lines 1-3, an element from the input polynomial set $F$ is
chosen such that its leading monomial has no proper divisors (in the conventional,
non-involutive sense)
among the remaining elements in $\lm(F)$. Similar choice of an element from $Q$ is
done in line
7. Before insertion of a new nonzero element $h$ into $G$ done in line 15,
all elements $g\in G$ such that $\lm(g)$ is a proper multiple of $\lm(h)$ are moved from
$G$ to $Q$ in line 13.

Such a choice in lines 1 and 7 together with the displacement step in line 13 provide
correctness of the algorithm. If it terminates, then the output basis
obviously satisfies the conditions~(\ref{inv_bas_1}) since $Q$ becomes the empty set
by the {\bf while}-condition in line 18. Minimality of the output basis can be
proved~\cite{G02-1} by exactly the same arguments as used in the proof of minimality for
the algorithm in paper~\cite{GB98b}.

As to termination of algorithm {\bf InvolutiveBasis I}, it follows from
Noetherianity of
division ${\cL}$. Indeed, there can be only finitely many cases~\cite{GB98a} when a nonzero
element $h$ obtained in line 9 and the polynomial $p$ selected in line 7 is ${\cL}-$head reducible.
Between such events, the set $\lm{G}$ of leading monomials
in $G$ can only be completed in the main loop 4-18 by finitely many monomials which are not
non-multiplicative prolongations of monomials in $\lm{G}$. All other completions are performed
just by these prolongations. By the Noetherianity of ${\cL}$, the last completion must
terminate in finitely many steps~\cite{GB98a}.

\vskip 0.3cm
\noindent
\begin{algorithm}{InvolutiveBasis I ($F,\prec,{\cL}$)\label{InvolutiveBasis}}
\begin{algorithmic}[1]
\INPUT $F\in \R\setminus \{0\}$, a finite set;\ $\prec$, an order;\
   ${\cL}$, an involutive division
\OUTPUT $G$, a minimal involutive basis of $\Id(F)$
\STATE {\bf choose} $f\in F$ without proper divisors of $\lm(f)$ in
 $\lm(F)\setminus \{\lm(f)\}$
\STATE $G:=\{f\}$
\STATE $Q:=F\setminus G$
\DOWHILE
  \STATE $h:=0$
  \WHILE{$Q\neq \emptyset$\ and $h=0$}
    \STATE {\bf choose} $p\in Q$ without proper divisors of $\lm(p)$ in
 $\lm(Q)\setminus \{\lm(p)\}$
    \STATE $Q:=Q\setminus \{p\}$
    \STATE $h:={\bf NormalForm}(p,G,\prec,{\cL})$
  \ENDWHILE
  \IF{$h\neq 0$}
    \FORALL{$\{ g\in G \mid \lm(h)\sqsubset \lm(g)\}$}
      \STATE $Q:=Q\cup \{g\}$; \ $G:=G\setminus \{g\}$
    \ENDFOR
    \STATE $G:=G\cup \{ h \}$
    \STATE $Q:=Q\cup \{\,g\cdot x \mid g\in G,\,
                                    x\in NM_{\cL}(\lm(g),\lm(G))\,\}$
 \ENDIF
\ENDDO{$Q \neq \emptyset$}
\RETURN $G$
\end{algorithmic}
\end{algorithm}

\vskip 0.3cm
\noindent
It should be noted that the above algorithm is distinguished from the algorithm
in~\cite{GBY01} by the conditions of the choice and of the displacement used in lines 1, 3
and 7, respectively. In paper~\cite{GBY01} the polynomials with the smallest leading
monomials with respect to the order $\succ$ were chosen, and those with higher leading
monomials then $\lm(h)$ were displaced. The choice made in~\cite{GBY01} is apparently more
restrictive than that in the above algorithm {\bf InvolutiveBasis I}, and with more
replacements than those done in line 13 of the algorithm.

\vskip 0.3cm
\noindent
\begin{algorithm}{NormalForm$(p,G,\prec,{\cL})$}
\begin{algorithmic}[1]
\INPUT $p\in \R\setminus \{0\}$, a polynomial; $G\subset \R\setminus \{0\}$, a finite set;
  \\ \hspace*{0.6cm} $\prec$, an order;\ ${\cL}$, an involutive division
\OUTPUT $h=NF_{\cL}(p,G)$, the ${\cL}-$normal form of $p$ modulo $G$
      \STATE $h:=p$
      \WHILE{$h\neq 0$ {\bf and} $h$ has a term $t$ ${\cL}-$reducible modulo $\lm(G)$}
            \STATE {\bf take} $g\in G$ such that $\lm(g)\mid_{\cL} t$
            \STATE $h:=h - g\cdot t/\lt(g)$
      \ENDWHILE
  \RETURN $h$
\end{algorithmic}
\end{algorithm}

\vskip 0.3cm
\noindent
Subalgorithm {\bf NormalForm} invoked in line 9 of the {\bf InvolutiveBasis I} algorithm
computes ${\cL}-$normal form in the full
accordance with Definition 2.6 specified
for ${\cL}-$division. Its termination immediately follows from the fact that
${\cL}-$reductions form a subset of the conventional reductions, and the last reduction
sequence is always finite~\cite{BB65,Buch85}.

${\cL}-$reducibility of polynomial $h$ is checked in line 2 of algorithm {\bf NormalForm}.
The check
consists in search for an ${\cL}-$divisor of $\lm(h)$ among elements in
$\lm(G)$. Polynomial set $G$, as constructed in the course of the main algorithm
{\bf InvolutiveBasis I}, is ${\cL}-$autoreduced
at every step of the completion procedure~\cite{GB98b,G02-1}. Therefore, by property 2
in Definition 3.1, $\lm(h)$ may have at most one ${\cL}-$divisor. Hence,
polynomial $g$ satisfying the condition taken in line 3 is unique in $G$
as an ${\cL}-$reductor of the term $t$ in $h$.

Now we illustrate algorithm {\bf InvolutiveBasis I} by Example 3.4. In Table 2 we show
the intermediate values of sets $G$ in the 2nd column, Janet non-multiplicative variables
$NM_J$ for elements in $G$ in the 3rd column and elements in $Q$ in the 4th column. Rows
of the table contain these values obtained at the initialization and after every
iteration of the main loop 4-18. In this case in lines 3 and 7 we selected the lexicographically
smallest elements.

\noindent
\begin{table}[h]
\caption{Computation of Janet basis for Example 3.4}
\begin{center}
\begin{tabular} {|c|c|c|c|} \hline\hline
Steps of  & \multicolumn{3}{c|}{Sets $G$ and $Q$}
\\ \cline{2-4}
 algorithm & elements in $G$ & $NM_J$ & $Q$      \\ \hline
 initialization & $xy^2-1$ & $-$ & $\{x^2y-1\}$ \\ \hline
 iteration  & $x^2y-1$ & $-$ &   \\
            & $xy^2-1$ & $x$ & $\{x^2y^2-x\}$ \\ \cline{2-4}
            & $x-y$    & $-$ & $\{xy^2-1,x^2y-1\}$   \\ \cline{2-4}
            & $x-y$    & $-$ &  \\
            & $y^3-1$  & $x$ & $\{x^2y-1,xy^3-x\}$   \\ \cline{2-4}
            & $x-y$    & $-$ &  \\
            & $y^3-1$  & $x$ & $\{\ \}$   \\ \hline \hline
\end{tabular}
\end{center}
\end{table}

\noindent
The grave practical disadvantage of the presented algorithm is that it treats useless
repeated prolongations and does not use any criteria to avoid unnecessary reductions.
Below we describe an improved version of involutive algorithm where unnecessary
repeated prolongations are avoided and where the involutive analogues of Buchberger's
criteria are enabled. For these purposes we need the next definition.

\vskip 0.3cm
\noindent
{\bf Definition 4.1.}\ \ An {\em ancestor} of a polynomial $f\in F\subset \R\setminus \{0\}$ is a
polynomial $g\in F$ of the smallest $\deg(\lm(g))$ among those satisfying
$f=g\cdot u$ modulo $\Id(F\setminus \{f\})$
with $u\in \M$. If $\deg(\lm(g))<\deg(\lm(f))$ ($u\neq 1$) the ancestor $g$ of $f$ is
called {\em proper}.
\noindent
If an intermediate polynomial $h$ that arose in the course of a completion algorithm
has a proper ancestor $g$, then $h$ has been obtained from $g$ via a sequence
of ${\cL}$-head irreducible non-multiplicative prolongations. For the ancestor $g$
itself the equality $\lm(\anc(g))=\lm(g)$ holds.

Let now every element
$f\in F$ in the intermediate set of polynomials be endowed (cf.~\cite{GBY01})
with the triple structure
$$
p=\{f,\, g,\, vars\}
$$
where
$$
\begin{array}{lcl}
\pol(p)&=&f\ \ \mbox{is the polynomial}\ f\ \mbox{itself},\\
\anc(p)&=&g\ \ \mbox{is a polynomial ancestor of}\ f\
 \mbox{in}\ F,\\
\nmp(p)&=&vars\ \ \mbox{is a (possibly empty) subset of variables}.
\end{array}
$$
The set $vars$ associated with polynomial $f$ accumulates those non-multiplicative
variables of $f$ have been already used in the algorithm for construction of
non-multiplicative prolongations. It keeps information on non-multiplicative
prolongations of polynomial $f$ that have been already examined
in the course of completion and serves to avoid useless repeated prolongations.

Knowledge of an ancestor of $f$ in $F$ helps to avoid some unnecessary reductions by
applying the involutive analogues of Buchberger's criteria described below.
An improved version of algorithm {\bf InvolutiveBasis I} named {\bf InvolutiveBasis II} is
given as follows.

\vskip 0.3cm
\noindent
\begin{algorithm}{InvolutiveBasis II($F,\prec,{\cL}$)\label{InvBasImpr}}
\begin{algorithmic}[1]
\INPUT $F\in \R\setminus \{0\}$, a finite set; $\prec$, an order; ${\cL}$, an involutive division \\
\OUTPUT $G$, a minimal involutive  basis of $\Id(F)$ or a reduced \Gr basis
\STATE {\bf choose} $f\in F$ without proper divisors of $\lm(f)$ in
 $\lm(F)\setminus \{\lm(f)\}$
\STATE $T:=\{f,f,\emptyset\}$
\STATE $Q:=\{\{q,q,\emptyset\} \mid q\in F\setminus \{f\}\}$
\STATE $Q:=${\bf HeadReduce}$(Q,T,\prec,{\cL})$
  \WHILE{$Q \neq \emptyset$}
    \STATE {\bf choose} $p \in Q$ without proper divisors of $\lm(\pol(f))$ in \\
    \hspace*{1.2cm} $\lm(\pol(Q))\setminus \{\lm(\pol(f))\}$
    \STATE $Q:=Q\setminus \{p\}$
    \IF{$\pol(p) = \anc(p)$}
          \FORALL{$\{q \in T \mid \lm(\pol(q))\sqsupset \lm(\pol(p))\}$}
          \STATE $Q:=Q \cup \{q\}$; \hspace*{0.4cm} $T:=T \setminus \{q\}$
          \ENDFOR
    \ENDIF
    \STATE $h:=${\bf{TailNormalForm}}($p,T,\prec,{\cL})$
    \STATE $T:=T \cup \{h,\anc(p),\nmp(p)\}$
    \FORALL{$q\in T$ {\bf and} $x\in NM_{\cL}(q,T)\setminus \nmp(q)$}
      \STATE $Q:=Q \cup \{\{\pol(q)\cdot x,\anc(q),\emptyset\}\}$
      \STATE $\nmp(q):=\nmp(q)\cap NM_{\cL}(q,T)\cup \{x\}$
    \ENDFOR
    \STATE $Q:=${\bf HeadReduce}$(Q,T,\prec,{\cL})$
  \ENDWHILE
  \RETURN $\{\pol(f)\mid f\in T\}$ or $\{\pol(f)\mid f\in T\mid f=\anc(f)\}$
\end{algorithmic}
\end{algorithm}

\vskip 0.3cm
\noindent
Here and in the below algorithms, where no confusion can arise,
we simply refer to the triple set $T$ as the second argument in $NM_{\cL}$, $NF_{\cL}$,
and $HNF_{\cL}$
instead of the polynomial set $\{\,g=\pol(t) \mid t\in T\,\}$. Sometimes we also refer to the
triple $p$ instead of $\pol(p)$. Besides, when we speak on reduction of triple set $Q$ modulo
triple set $T$ we mean reduction of the polynomial set $\{\,f=\pol(q) \mid q\in Q\,\}$
modulo $\{\,g=\pol(t) \mid t\in T\,\}$.

Apart from providing subsets $Q$ and $T$ of the intermediate basis with
the triple structure, the improved version contains extra lines 4 and 19. Here ${\cL}-$
head reduction is done for the basis elements in $Q$ modulo those in $T$. Then
the remaining tail reduction is performed in line 13 to obtain the (full)
${\cL}$-involutive normal form. Furthermore, unlike the previous algorithm and
due to presence of the third elements in triples, the set $Q$ is enlarged in line 16 only
with those non-multiplicative prolongations which have not been examined yet.

In doing so, a new prolongation of a polynomial in $T$ is inserted into $Q$ with the
ancestor of the polynomial. In the next line the selected non-multiplicative variable
$x$ is added to the set of non-multiplicative variables already used. The intersection in line 17
takes into account that some of these variables may turn into multiplicative
through the contraction of $T$ in line 10 and by virtue of the relation 3
in Definition 3.1.

The subalgorithm {\bf HeadReduce} invoked in lines
4 and 19 of the main algorithm returns set $Q$ which, if nonempty, contains the part of
intermediate basis ${\cL}-$head reduced modulo $T$.

\vskip 0.3cm
\noindent
\begin{algorithm}{HeadReduce$(Q,T,\prec,{\cL})$}
\begin{algorithmic}[1]
\INPUT $Q$ and $T$, sets of triples; $\prec$, an order;\ ${\cL}$, an involutive division
\OUTPUT ${\cL}-$head reduced set $Q$ modulo $T$
 \STATE $S:=Q$
 \STATE $Q:=\emptyset$
 \WHILE {$S \neq \emptyset$}
   \STATE {\bf choose} $p\in S$
   \STATE $S:=S\setminus \{p\}$
   \STATE $h:=${\bf{HeadNormalForm}}$(p,T,{\cL})$
   \IF{$h\neq 0$}
       \IF{$\lm(\pol(p))\neq \lm(h)$}
         \STATE $Q:=Q\cup \{h,h,\emptyset\}$
       \ELSE
         \STATE $Q:=Q\cup \{p\}$
       \ENDIF
   \ELSE
      \IF{$\lm(\pol(p))=\lm(\anc(p))$}
         \FORALL{$\{ q\in S \mid \anc(q)=\pol(p)\}$}
          \STATE $S:=S\setminus \{q\}$
         \ENDFOR
      \ENDIF
   \ENDIF
\ENDWHILE \RETURN $Q$
\end{algorithmic}
\end{algorithm}

\vskip 0.3cm
\noindent
Its own subalgorithm {\bf HeadNormalForm}
invoked in line 6 just computes the ${\cL}-$head normal form $h$ of the polynomial $\pol(p)$
in the input set $Q$ assigned to the set $S$ at the initialization step (lines 1-2). If $h\neq 0$
when $\lm(\pol(p))$ is ${\cL}-$reducible what is verified in line 8, then $\lm(h)$ does not belong
to the initial ideal generated by $\{\,\lm(\pol(f)) \mid f\in Q\cup T\,\}$~\cite{GB98b}.
In this case the triple $\{h,h,\emptyset\}$ for $h$ is inserted
(line 9) into the output set $Q$. Otherwise, the output set $Q$ retains the triple $p$ as it
is in the input.

In the case when $h=0$, whereas $\pol(p)$ has no proper ancestors what is verified in line 14,
all the descendant triples for $p$, if any, are deleted from $S$ in line 16. Such
descendants cannot occur in $T$ owing to the choice conditions in lines 1, 6  and to
the displacement condition of line 9 in the main algorithm {\bf InvolutiveBasis II}. Step 14-18
serves for the memory saving and can be ignored if the memory restrictions are not very
critical for a given problem. In this case all those descendants will be casted away
by the criteria checked in the below algorithm {\bf HeadNormalForm}.

The next algorithm {\bf HeadNormalForm}, after initialization in lines 1-2,
starts with verification in line 3 of ${\cL}-$head reducibility of the input
polynomial $h$ modulo polynomial set $G:=\{\,\pol(g)\mid g\in T\,\}$.
This verification consists in search for ${\cL}-$divisor of $\lm(h)$ in $\lm(G)$.
If there is no such divisor the algorithm returns $h$ in line 4. Otherwise,
in the course of the search, the polynomial $g\in G$
is found such that $\lm(g)\mid_{\cL} lm(h)$.

\vskip 0.3cm
\noindent
\begin{algorithm}{HeadNormalForm($p,T,\prec,{\cL})$}
\begin{algorithmic}[1]
\INPUT $T$, a set of triples; $p$, a triple such that $\pol(p)=\pol(f)\cdot x$, \\
 \hspace*{0.4cm} $f\in T$, $x\in NM_{\cal{L}}(f,T)$;
 $\prec$, an order; ${\cL}$, an involutive division
\OUTPUT $h=HNF_{\cL}(p,T)$, the ${\cL}-$head normal form
   of $\pol(p)$ modulo $T$
   \STATE $h:=\pol(p)$
   \STATE $G:=\{\,\pol(g)\mid g\in T\,\}$
    \IF{$\lm(h)$ is ${\cL}$-irreducible modulo $G$}
       \RETURN $h$
    \ELSE
       \STATE {\bf take} $g\in T$ such that $\lm(\pol(g))\mid_{\cL} \lm(h)$
       \IF{$\lm(h) \neq \lm(\anc(p))$}
          \IF{{\em Criteria}$(p,g)$}
            \RETURN $0$
          \ENDIF
       \ELSE
         \WHILE{$h\neq 0$ {\bf and} $\lm(h)$ is ${\cL}-$reducible modulo $G$}
            \STATE {\bf take} $q$ in $T$ such that $\lm(q)\mid_{\cL} \lm(h)$
            \STATE $h:=h - q\cdot \lt(h)/\lt(q)$
         \ENDWHILE
       \ENDIF
    \ENDIF \ \ {\bf return $h$}
\end{algorithmic}
\end{algorithm}

\vskip 0.3cm
\noindent
If the set $G$ is ${\cL}-$autoreduced, and this just takes place when algorithm {\bf HeadNormalForm} is
invoked in line 6 of algorithm {\bf HeadReduce}, then there is the only one such $g$ in
$G$. Similarly, the check whether the further ${\cL}-$reducibility takes place is done
in line 12, and if that is the case, then the corresponding unique reductor is taken in
line 13.

For the ${\cL}$-head reducible input polynomial $\pol(p)$ what is checked
in line 3, the following four criteria are verified in line 8
\vskip 0.2cm
\centerline{\em Criteria$(p,g)$=$C_1(p,g) \vee C_2(p,g) \vee  C_3(p,g) \vee C_4(p,g)$,}
\vskip 0.2cm
\noindent
where
\begin{itemize}

\item[] $C_1(p,g)$ is true $\Longleftrightarrow$ $\lm(\anc(p))\cdot \lm(\anc(g))
 = \lm(\pol(p))$,

\item[] $C_2(p,g)$ is true $\Longleftrightarrow$
$\lcm(\lm(\anc(p)), \lm(\anc(g))) \sqsubset \lm(\pol(p)) $,

\item[] $C_3(p,g)$ is true $\Longleftrightarrow$ $\exists$ $t\in T$ such that \\
\hspace*{1.0cm}$\lcm(\lm(\pol(t)),\lm(\anc(p))) \sqsubset \lcm(\lm(\anc(p)), \lm(\anc(g)))$
$\wedge$ \\
\hspace*{1.0cm}$\lcm(\lm(\pol(t)),\lm(\anc(g))) \sqsubset \lcm(\lm(\anc(p)), \lm(\anc(g)))$,

\item[] $C_4(p,g)$ is true $\Longleftrightarrow$ $\exists$ $t\in T$ $\wedge$
$y\in NM_{\cL}(t,T)$ with $\lm(\pol(t))\cdot y=\lm(\pol(p))$,  \\
\hspace*{1.0cm}$\lcm(\lm(\anc(p)),\lm(\anc(t))) \sqsubset \lm(\pol(p))$
$\wedge$ $\idx(t,T)<\idx(f,T)$,
\end{itemize}
where $\idx(t,T)$ enumerates position of triple $t$ in set $T$.

Criterion $C_1$ is Buchberger's co-prime criterion~\cite{Buch85} in
its involutive form~\cite{GBY01}. It is easy to see that under the
condition $\lm(\pol(g))\mid_{\cL} \lm(h)$ of line 6
$$\lm(\anc(f))\cdot \lm(\anc(g)) = \lm(\pol(f))\Longleftrightarrow \lm(\anc(f))\cdot \lm(\anc(g)) \mid \lm(\pol(f)),$$
and in~\cite{GBY01} the right-hand side form of the criterion was given.
\noindent
Criterion $C_2$ was derived in~\cite{GB98a} as a consequence of Buchberger's chain
criterion~\cite{Buch85}. Criteria $C_3$ and $C_4$ derived in
paper~\cite{AH02} complement criterion $C_2$ to the full equivalence to the
chain criterion. Thus, the four criteria in the aggregate are
equivalent~\cite{AH02} to Buchberger's criteria adapted to the involutive completion
procedure.

The last subalgorithm {\bf TailNormalForm} completes the ${\cL}-$head reduction
by performing the involutive tail reduction and is invoked in line 13 of the main
algorithm {\bf InvolutiveBasis II}. In that way, it returns the ${\cL}-$normal form of
the input ${\cL}-$head reduced polynomial as given by Definition 2.6.

\vskip 0.3cm
\noindent
\begin{algorithm}{TailNormalForm$(p,T,\prec,{\cL})$}
\begin{algorithmic}[1]
\INPUT $p$, a triple $p$ such that
   $\pol(p)=HNF_{\cL}(p,T)$; \\ \hspace*{0.6cm} $T$, a set of triples; $\prec$, an order;\ ${\cL}$, an involutive division
\OUTPUT $h=NF_{\cL}(p,T)$, the ${\cL}-$normal form of $\pol(p)$ modulo $T$
    \STATE $G:=\{\,\pol(g)\mid g\in T\,\}$
    \STATE $h:=\pol(p)$
      \WHILE{$h$ has a term $t$ ${\cL}-$reducible modulo $G$}
            \STATE {\bf take} $g\in G$ such that $\lm(g)\mid_{\cL} t$
            \STATE $h:=h - g\cdot t/\lt(g)$
      \ENDWHILE
  \RETURN $h$
\end{algorithmic}
\end{algorithm}

\vskip 0.3cm
\noindent
The main algorithm can output either involutive or reduced \Gr basis (or both) depending on
the instruction used in line 21. It is easy to verify that the main algorithm
together with its subalgorithms ensures that every element in the output involutive basis
has one and only one ancestor. This ancestor is apparently irreducible, in the \Gr sense, by
other elements in the basis. Thereby, those elements in the involutive basis that have no
proper ancestors constitute the reduced \Gr basis.

Though involutive bases are usually redundant as \Gr ones, this specific redundancy,
can provide more accessibility to information on polynomial ideals and modules~\cite{S-2,Apel03}.
For example, unlike reduced \Gr bases, involutive bases give explicit simple formulae for
both the Hilbert function~\cite{Apel} and Hilbert polynomial~\cite{G99,GBY01} of a
polynomial ideal ${\cal{I}}$. If $G$ is an involutive basis of ${\cal{I}}$, then
the (affine) Hilbert function $HF_{{\cal{I}}}(s)$ and the Hilbert polynomial
$HP_{{\cal{I}}}(s)$ are expressed in terms of binomial coefficients as follows.
\begin{eqnarray*}
HF_{{\cal{I}}}(s)&=&\left(\begin{array}{c} n+s \\ s\end{array}\right) -
 \sum_{i=0}^s \sum_{u\in \lm(G)} \left(\begin{array}{c} i-deg(u)+\mu(u)-1 \\
\mu(u)-1 \end{array}\right)\,, \\[1.0cm]
HP_{{\cal{I}}}(s)&=&\left(\begin{array}{c} n+s \\ s\end{array}\right) -
 \sum_{u\in \lm(G)} \left(\begin{array}{c} s-deg(u)+\mu(u) \\
\mu(u) \end{array}\right)\,,
\end{eqnarray*}
where $\mu(u)$ is the number of multiplicative variables for $u$.
The first term in the right hand sides of these formulae is the
total number of monomials in $\M$ of degree $\leq s$. The sum in
the expression for $HF_{{\cal{I}}}(s)$ counts~\cite{Apel} the
number of monomials of degree $\leq s$ in the ${\cL}-$cone (see
Definition 3.2)\ \ ${\cal{C}}_{\cL}(\lm(G))$. This number
coincides with the number of such monomials in the cone
${\cal{C}}(\lm(G))$ by the completeness condition
${\cal{C}}_{\cL}(\lm(G))={\cal{C}}(\lm(G))$ for $\lm(G)$ that
follows from Definition 3.5.

There is a number of other useful applications of involutive bases in commutative and noncommutative
algebra (see~\cite{S-1,S-2,Apel03} and references therein).

We conclude this section by observing that if one applies algorithm {\bf InvolutiveBasis II} to
construction Janet basis for the polynomial set in Example 3.4, then the intermediate
polynomial data coincide with those in Table 2, where $G=\{\, \pol(p)\mid p\in T\,\}$ denotes
the $T-$part of the intermediate basis and $Q$ stands for its $Q-$part. The only distinction from
functioning algorithm {\bf InvolutiveBasis I} for this particular example is that
prolongation $xy^3-x$ reduces to zero by criterion $C_1$.

\section{Efficiency Issues}

Computer experiments presented on the Web cite {\tt http://invo.jinr.ru}
reveal superiority (for most of benchmarks) of the C/C++ code~\cite{GBY01} implementing
the {\bf InvolutiveBasis II} algorithm over the best present-day implementations of
Buchnerger's algorithm.

In this section we consider some efficiency aspects of  algorithm
{\bf InvolutiveBasis II} as compared with the Buchberger's algorithm~\cite{Buch85}
when one treats all possible $S-$polynomials as critical pairs, under application of the both
Buchberger criteria, and performs reductions by means of the conventional monomial division.

\subsection{Automatic avoidance of some useless critical pairs}

In the involutive case, whenever the leading monomial of the intermediate polynomial $f$,
under consideration in the course of a completion algorithm, is ${\cL}-$ reducible
modulo another intermediate polynomial $g$, the corresponding $S-$polynomial whose
${\cL}-$normal form is to be computed has the structure $S(f,g)=f-g\, u$. Monomial
$u\in \M$ contains only multiplicative variables for $g$. In terms of the ancestors
$\tilde{f}$ and $\tilde{g}$ for $f$ and $g$, respectively, this
$S-$polynomial corresponds to the conventional $S-$polynomial
\begin{equation}
S(\tilde{f},\tilde{g})=\tilde{f}\,\frac{\lm(f)}{\lm(\tilde{f})}-
\tilde{g}\,\frac{\lm(g)}{\lm(\tilde{g})}\,u\,. \label{inv_spol}
\end{equation}
If ancestor $\tilde{f}$ is proper, then $f$ has been obtained from $\tilde{f}$ by a sequence of
${\cL}-$head irreducible non-multiplicative prolongations. Insertion of these prolongations into the
intermediate basis can only be accompanied with tail reductions.

Thus, in the involutive procedure only $S-$polynomials (critical pais) of form (\ref{inv_spol})
must be processed by computing their involutive normal form. In this way some useless
$S-$polynomials are automatically ignored in the course of involutive algorithms.

As a simple example, consider the polynomial set presented in Table 3.
The first column contains the polynomials. Their Janet non-multiplicative
variables for $x\succ y\succ z$ and the corresponding prolongations shown in the second and
third column, respectively. The forth column contains $S-$polynomials to be examined in a completion
procedure. We see that $S-$polynomial $S(p_2,p_3)$ is not considered at all. To avoid
useless reduction of $S(p_2,p_3)$, Buchberger's algorithm needs the chain criterion.

\begin{table}
\caption{Avoidance of useless critical pair in involutive
completion}
\begin{center}
\begin{tabular} {|c|c|c|c|} \hline \hline
Polynomial & $NM_J$ & Prolongation & $S-$polynomial       \\ \hline
 $p_1=xy-1$ & $-$ & $-$      & $-$   \\
 $p_2=xz-1$ & $y$ & $y\,p_2$ & $S(p_2,p_1)=y\,p_2 - z\,p_1$\\
 $p_3=yz-1$ & $x$ & $x\,p_3$ & $S(p_3,p_1)=x\,p_3 - z\,p_1$\\ \hline \hline
\end{tabular}
\end{center}
\end{table}

\subsection{Weakened role of criteria}

The fact that a number of useless $S-$polynomials is automatically avoided in the involutive
approach contributes in its much weaker dependence on the use of criteria as compared with
Buchberger's algorithm. Already the first implementation of the involutive completion
procedure for zero-dimen\-sional ideals and Pommaret division done by Zharkov and Blinkov
in Reduce~\cite{ZB} revealed higher computation efficiency of this procedure in comparison
with the Reduce implementation of Buchberger's algorithm. And at that time it was unexpected,
since no criteria were used in~\cite{ZB} whereas Buchberger's criteria had been
implemented in the Groebner package of Reduce. It is well-known that without the criteria
Buchberger's algorithm becomes impractical even for rather small problems.

It should be also said that our implementation in~\cite{GBY01} of algorithm
{\bf InvolutiveBasis II} for Janet division, as presented on the Web site {\tt http://invo.jinr.ru},
exploited only criteria $C_1$ and $C_2$ (Sect.4). The effect of other two criteria $C_3$ and $C_4$
that are due to Apel and Hemmecke~\cite{AH02} (see also~\cite{HemThesis}) for our algorithm is now under
experimental investigation.

\noindent
\begin{table}[h]
\caption{Influence of criteria}
\begin{center}
\begin{tabular} {|c|c|c|c|c|r|r|r|r|} \hline\hline
 & \multicolumn{4}{c|}{Applicability} &\multicolumn{4}{c|}{Timing (sec.)}
 \\ \cline{2-9}
Example & $C_1$ & $C_2$ & $C_3$ & $C_4$ & Without & $C_1 \div C_2$  & $C_1 \div C_3$ & $C_1 \div C_4$ \\ \hline
Cyclic6  & 98  &  2     &   4   &  --     &  0.14  &  0.09  &   0.11  &  0.12 \\
Cyclic7  & 698 & 190    &  22   &  --    & 73.01  &  46.54 &  47.15  &  58.72 \\
Katsura8 & 173 & 1      &   1   &  --    & 29.54  &  22.44 &  22.25  &  27.48 \\
Katsura9 & 344 & --     &   1   &  --    & 358.34 & 282.25 & 278.25  &  337.52\\
Cohn3     & --  & 114    & 169  &  7     & 736.15 &  81.87 &  55.10  &  76.72 \\
Assur44   & 89  &  60    & 171  &  3    & 10.76  &  10.13 &   8.83  &  10.35 \\
Reimer6  & 63  & 235    & 179   &  12    & 30.65  &  17.06 &   7.46  &   9.69 \\ \hline \hline
\end{tabular}
\end{center}
\end{table}

\noindent
Table 4 contains data of computer experiments with the
criteria for some benchmarks from~\cite{BM,Verschelde} and for computation over the ring
of integers. Columns 2-5 in Table 4 show how many times the criteria are applied in our algorithm
for the degree-reverse-lexicographical order. In so doing, criterion $C_1$ is checked first,
then sequentially $C_2$, $C_3$ and $C_4$~\footnote{The numbers in columns 2-5 of Table 4 for Cyclic6
widely differ from those in~\cite{HemThesis} where another involutive algorithm was
experimentally analyzed.}.

The last four columns of Table 4 illustrate the effect of the criteria for the C code~\cite{GBY01}
running on an Opteron-242 computer under Gentoo Linux. Comparison of the timings with and without
criteria shows that the use of criteria is not so critical for algorithm {\bf InvolutiveBasis II}
as for Buchberger's algorithm. For the latter problems of the Cyclic7 size and larger become
intractable if no criteria applied.

Note, that the check of criterion $C_4$ is more expensive from the computational
point of view~\cite{HemThesis}.  That is why for all the examples in the table
thsi criterion is adversely affected on the timings as Table 4 indicates.

In any event, however, criteria $C_3$ and $C_4$ comprise qualitative improvement of
the involutive algorithm and their quantitative effect has to be further studied.

\subsection{Smooth growth of intermediate coefficients}

We performed experimental investigation of intermediate coefficient growth for
polynomial systems with integer coefficients as they are completed to involution by
algorithm {\bf InvolutiveBasis II} for Janet division. Our observation is that this
growth is much more smooth than in the case of Buchberger's algorithm. Swell of
intermediate coefficients is a well-known difficulty of this algorithm.
Even in the case when coefficients in the initial polynomial set and in the \Gr basis
are small, intermediate coefficients can be huge and lead to a dramatic slowing down
of computation or to running out of a computer memory. The following example taken
from~\cite{Arnold} nicely illustrates this behavior of Buchberger's algorithm

\vskip 0.3cm
\noindent
{\bf Example 5.1.}\ Consider ideal ${\cal{I}}=Id(F)$ in $\Q[x,y,z]$ generated by the
polynomial set:
\begin{equation}
F=\left\{
\begin{array}{l}
 8\,x^2y^2 + 5x\,y^3 + 3x^3z + x^2y\,z, \\[0.1cm]
 x^5 + 2\,y^3z^2 + 13\,y^2\,z^3 + 5\,y\,z^4, \\[0.1cm]
 8\,x^3 + 12\,y^3 + x\,z^2 + 3,\\[0.1cm]
 7\,x^2y^4 + 18\,x\,y^3z^2 + y^3z^3.
\end{array}
\right. \label{Arnold}
\end{equation}
Its \Gr basis for the degree-reverse-lexicographical order with
$x\succ y\succ z$ is small $G=\{x,4\,y^3+1,z^2\}$ whereas in the
course of Buchberger's algorithm, as it implemented in Macaulay 2,
there arise intermediate coefficients with about 80,000
digits~\cite{Arnold}. As to algorithm {\bf InvolutiveBasis II}, it
outputs \Gr basis $G$ or Janet basis $G\cup \{y\,z^2,y^2z^2\}$,
depending on the instruction in line 21, with not more than 400
digits in the intermediate coefficients. For comparison with the
timings for this example given in~\cite{Arnold}, we run our
code~\cite{GBY01} on a 500 Mhz PC with 256 Mb RAM and computed the
\Gr basis in 0.1 seconds.

\vskip 0.3cm
\noindent
For most of benchmarks from~\cite{BM,Verschelde} and the degree-reverse-lexicographic order
the intermediate coefficients arising in the course of the involutive algorithm
grow up only in several times in comparison with coefficients in the output basis. For the
same collection of examples as in Table.4, the maximal lengths of the input, intermediate
and output coefficients, measured by number of 64 bit words occupied, are accumulated in
columns
2-4 of Table 5. Here we applied the first three criteria without the forth one.
The last column shows the ratio of the length in the third and forth columns called
``swell factor''.

\noindent
\begin{table}[h]
\caption{Coefficient size in 64 bit words}
\begin{center}
\begin{tabular} {|c|c|c|c|c|} \hline\hline
Example   & Input  & Intermediate  & Output &  Swell factor \\ \hline
Cyclic6  &  1  &  3     &   1   &  3.00   \\
Cyclic7  &  1  & 11     &   5   &  2.20   \\
Katsura8 &  1  &  5     &   4   &  1.25   \\
Katsura9 &  1  &  8     &   6   &  1.33  \\
Cohn3     &  1  & 168    &  19   &  8.84  \\
Assur44   &  1  &  93    &  19   &  4.89  \\
Reimer6  &  1  &  4     &   4   &  1.00  \\ \hline \hline
\end{tabular}
\end{center}
\end{table}

\noindent
Only for a few examples among those we have already analyzed, the swell factor reached
several dozens. For instance, for the ``f855'' benchmark it is 30.

To our opinion, there are several peculiarities
of the involutive completion procedure that may provide such a smooth behavior of the
intermediate arithmetics.
\begin{itemize}

\item Both selection of prolongations from $Q$ (and the critical pairs of form (\ref{inv_spol})
among them) for the reduction process and the reduction process itself are more
restrictive in the involutive approach than in Buchberger's algorithm. For the last
algorithm, as a result of numerous computer experiments, some strategies such as normal
strategy~\cite{Buch85} and mainly the ``sugar''~\cite{Traverso} we found as heuristically
``good'' restrictions for selection of critical pairs. In the above involutive algorithms
the leading term of a prolongation selected for the tail reduction and insertion into $T$
must not divide the leading terms of other prolongations. However, this restriction is
far yet from fixing the prolongation to be selected. For the degree-reverse-lexicographic order,
in our implementation~\cite{GBY01} we select that with the minimal total degree of the
leading monomial. By analogy with strategy in~\cite{Buch85}, this strategy can be considered
as normal.

In doing so, the ${\cL}-$head pre-reduction of $Q$ modulo $T$ done before selection of
a prolongation for insertion into $T$ is also an ingredient of the selection procedure. It
should be noted that in the {\bf HeadReduce} algorithm one can also do
the ${\cL}$-head reduction only for a part of elements in $Q$. Then selection in line 6
of the main algorithm must be done from the head reduced part. However, in the case of such
a partial ${\cL}$-head reduction of $Q$, the polynomial $\pol(p)$ in $p\in S$ chosen
in line 4 of algorithm {\bf HeadReduce} must not have proper divisors in
$\lm(\{\,\pol(q)\mid  q\in S\,\})\setminus \{\pol(p)\}$. The option to perform partial head
reduction gives rise to a certain freedom in the selection strategy. In our
implementation~\cite{GBY01} for the degree-reverse-lexicographic order we select
for the head reduction all the elements in $Q$ containing polynomials of the minimal
total degree.

As regards the reduction process itself, it is entirely fixed in the above involutive algorithms,
since there can be at most one elementary involutive reduction for any term. This is because of
uniqueness of an involutive divisor among the leading monomials of polynomial reductors in $T$
(see Sect.4).

\item The $T$-part of intermediate bases mostly contains some
extra polynomials, in comparison with intermediate bases in
Buchberger's algorithm that are redundant in the \Gr sense. These
extra polynomials which are non-multiplicative prolongations of
their ancestors in $T$, actively participate in the reduction
process and may also prevent intermediate coefficient swell as
already noticed in~\cite{Apel}.

\end{itemize}

\subsection{Fast search for involutive divisor}

As well as in any algorithm for constructing \Gr basis, the bulk of computing time for
examples large enough is expended for reductions. In the involutive reduction process the most
frequent operation is search for an involutive divisor for a given term among the
leading terms of the intermediate polynomial set. For the {\bf InvolutiveBasis II} algorithm
this operation is performed in lines 3 and 12 of subalgorithm  {\bf HeadNormalForm} and
in line 3 of subalgorithm {\bf TailNormalForm}. These two subalgorithms, especially
the former one,  are invoked (for examples large enough) enormous number of times from the main
algorithm and its subalgorithm {\bf HeadReduce}. Thus, an optimal search for ${\cL}$-divisor
is an essential ingredient of an efficient implementation of the involutive completion procedure.

We have already stated that in the above algorithms at every elementary reduction step there is
the only ${\cL}$-divisor. Therefore, a wanted search is such that when there is a divisor
it is located as quickly as possible. Otherwise, the search stops at the intermediate step
as early as possible and signals that there is no divisor.

For Janet division~\footnote{Pommaret division admits slightly different trees. For extension to other
divisions see~\cite{HemThesis} and \cite{Blinkov}.} we developed in~\cite{GBY01} a
special data structures called Janet
trees which allow to organize the wanted search. A Janet tree is a binary search tree
that takes proper account of properties of Janet division in Example 3.1 and whose
leaves contain (pointers at) monomials among which the searching is done.
Complexity bound for search in a Janet tree with the maximal total degree of the leaf
monomials $d$ in $n$ variables is $O(d+n)$. This is substantially lower than that for the
binary search in a sorted set of monomials~\cite{GBY01}.

We refer to~\cite{GBY01} for more details on Janet trees, and illustrate their
usefulness by an example.

\vskip 0.3cm
\noindent
{\bf Example 5.2.}\ Consider monomial set $U=\{x^3y,xz,y^2,yz,z^2\}$
and order $x\succ y\succ z$ on the variables. According to this order
which determines Janet separation of variables, assign number 1 to $x$,
number 2 to $y$ and number 3 to $z$. The below figure shows the structure of Janet tree
for set $U$.

\vskip 0.3cm
\noindent
The elements in $U$ are located in the leaves. Each edge is associated with
a certain variable, and each interior node has two integer indices. The root has
1 and 0 as its indices. For other interior nodes the first index
is the number of variable associated with the edge connecting the node with its
parent. The second index is the degree of the variable such that, if the node is not the root,
there exists a monomial in $U$ which contains this variable in the degree indicated.
In so doing, we consider the most compact form of the tree which is used by our
C/C++ code. The first index of the right interior child is
the number of the next variable. The left edge for the root is
associated with the first variable. And the left child of a node has always a
higher second index (degree in the current variable)
than its parent.

\vskip 0.3cm
\noindent
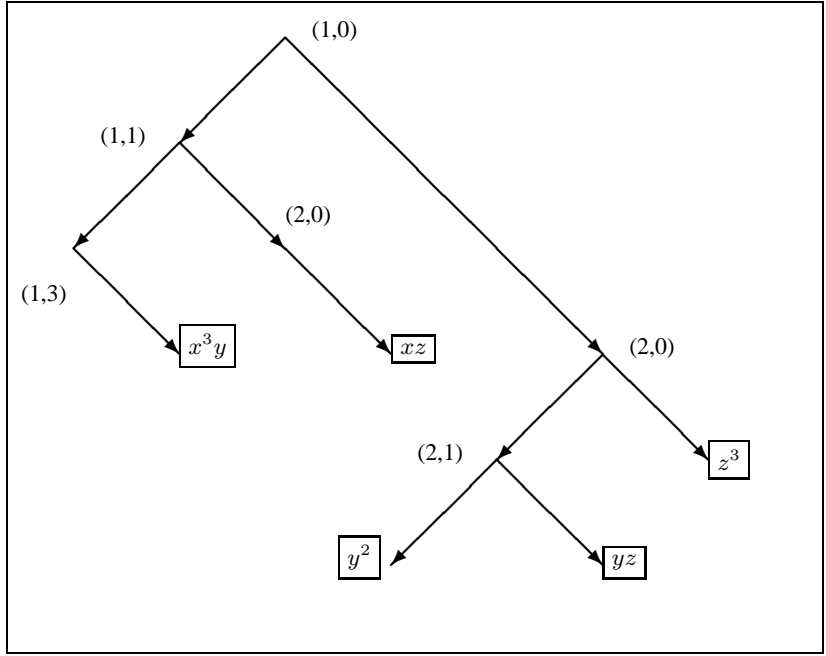
\begin{center}
{\small
\fbox{
\begin{picture}(300,240)(20,-40)
\thicklines \put(120,190){\vector(-1,-1){40}}
\put(80,150){\vector(-1,-1){40}}\put(40,110){\vector(1,-1){40}}
\put(80,70){\fbox{$x^3y$}} \put(80,150){\vector(1,-1){40}}
\put(120,110){\vector(1,-1){40}} \put(160,70){\fbox{$xz$}}

\put(120,190){\vector(1,-1){120}} \put(240,70){\vector(-1,-1){40}}
\put(200,30){\vector(1,-1){40}} \put(240,-10){\fbox{$yz$}}
\put(200,30){\vector(-1,-1){40}} \put(140,-10){\fbox{$y^2$}}

\put(240,70){\vector(1,-1){40}} \put(280,26){\fbox{$z^3$}}

\put(130,190){\small (1,0)} \put(50,150){\small (1,1)}
\put(20,90){\small (1,3)} \put(120,120){\small (2,0)}
\put(250,70){\small (2,0)} \put(170,30){\small (2,1)}
\end{picture}}
}
\end{center}

\vskip 0.3cm
\noindent
Given a monomial $u=x^iy^jz^k$ $(i,j,k\in \N{\geq 0})$, we start in our search for a
Janet divisor in $U$ from the root. If $i>0$ we select $x$ as the current variable
with number 1, and compare $i$ with the second index in the left child of the root
( 1 in our case ).

If $i=1$ we move to the child, and there is no further way to the left since the left
child of our current node has degree (second index) 3 in the current variable $x$.  Hence,
the leaf monomials of the lest subtree cannot divide $u$. Thus, we have to look
at the right child of our current node. If it exists, as in our case, we select the
second variable $y$ as the current one. Then we compare $j$ with the second index
(degree in $y$ and 0 in our case) of the right child. If $j\geq 0$ we move to the right child.
Since in our case the current node does not have left child, we change the current variable
to $z$ and compare $k$ with the degree in $z$ of the left child which is a leaf.
If $k\geq 1$ we go to the leaf and this means that the leaf monomial ($xz$ in our case)
is a Janet divisor we are looking for. Therefore for a monomial of the form
$xy^iz^k$ with $j\geq 0,k\geq 1$ Janet divisor is $xz$.

If $i=2$ we cannot move to the left, since for the current variable
$x$ the second index of the left child is $3$, and the leaf monomials in the left subtree
(in our case the only monomial $x^3y$) cannot contain Janet divisor of $u$.
Therefore, a monomial of the form $x^2y^jz^k$ has no Janet divisor in the tree.

If $i>2$ we go to the left child and then proceed with the next variable $y$, as
the current one, by comparing $j$ with the degree in $y$ of the right child. It is
leaf and has degree 1. It follows that $x^3y$ is a Janet divisor for $u$ with
$i\geq 3,j\geq 1,k\geq 0$ whereas for $i\geq 3,j=0,k\geq 0$ monomial $u$ has no
Janet divisor in $U$.

Similarly one can analyze the case $i=0$ when one has to go from the root to the
right.

This example illustrates the following general fact~\cite{GBY01}.
For a given monomial, a Janet tree provides the unique path in the tree which either
ends up with the leaf containing the sought divisor or breaks in a certain
interior node. The node is such that its left child, if any, has a higher degree in
the current variable and there is no way to the right since the right child has a
higher degree (in the next variable according to the order) than the given monomial.

\subsection{Parallelism}

Algorithm {\bf InvolutiveBasis II} admits a natural and efficient parallelism. Naturalness
is apparent from lines 4 and 19 of the main algorithm and from
the structure of subalgorithm {\bf HeadReduce}. The ${\cL}-$head reduction of polynomials
in $Q$, that is, the most time-consuming
part of the completion procedure can be done in parallel. As to efficiency,
in our recent paper~\cite{GY04} a slightly modified version of the algorithm oriented to
the multi-thread computation was experimentally studied on a two processor
Pentium III 700 Mhz computer running under Gentoo Linux. Some of experimental data
extracted from those obtained in~\cite{GY04} and related to the
benchmarks already used above, presented in Table 6.
\noindent
\begin{table}[h]
\caption{Timings (in seconds) and speedup due to parallelism}
\begin{center}
\begin{tabular}{|l|r|r|r|r|l|} \hline \hline
 Example & 1 Thread & 3 Threads & Speedup & $t_{1\,th}/t_{3\,th}$
\\ \hline
 Cyclic6     & 0.79  & 1.16          &  -0.37   &  0.68\ \ \ \ \\
 Cyclic7     & 386.89& 182.86        &  +294.03 &  2.12\ \ \ \ \\
 Katsura8    & 119.92 & 53.72        &  +66.20 &  2.23\ \ \ \ \\
 Katsura9    & 1356.37 & 587.82      &  +768.55&  2.31\ \ \ \ \\
 Cohn3       & 554.75& 222.69        &  +332.06&  2.49\ \ \ \ \\
 Assur44     & 73.93 & 31.34         &  +42.59 &  2.36\ \ \ \ \\
 Reimer6     & 88.99 & 52.56         &  +36.43  &  1.69\ \ \ \ \\
\hline \hline
\end{tabular}
\end{center}
\end{table}
\noindent
The second and third columns show the timings for the one- and three-thread modes of parallel
computation. The last two columns give absolute and relative speeding-up of the
three-thread run in comparison with one-thread. Each of the threads was doing the involutive
head reduction. We experimented with different numbers of treads, and it turned out that
the maximal speed-up is achieved just for three threads. This is in conformity with
the well-known observation that for multi-treading on a computer with SMP (Symmetric Multiprocessing)
architecture an optimal number of threads is the number of processors exceeded by one.

For such a small example as Cyclic6 the overheads of multi-threading surpass the
computational effect of the parallelization, and we obtain some slowing down. For the
other examples in Table 6 the speed-up is rather considerable. The relative speed-up
greater than the number of processors can be explained by a better selection
strategy dynamically realizable in the multi-thread mode. In this case some
head-reduced prolongations may come into play earlier than in the one-thread mode
and cause a faster chain of the successive reductions.

Since computational costs for Buchberger's algorithm are highly unstable
with respect to the selection strategy for $S-$polynomials, its experimental
parallelization (see, for example,~\cite{MTG}) does not reveal a reasonable
gain from the parallelism. The involutive algorithm, according to our two-processor
experiments, does not have this difficulty and rather stable to variation of the selection
strategy. Moreover, if one performs in parallel the ${\cL}-$head reduction,
then, as we observed in the three-thread run on the two-processor machine,
it may even optimize selection strategy. At least, for a degree compatible
term order, if the normal strategy is used, the parallelization apparently
helps to select a prolongation with the smallest total degree of the leading term.
As we observed experimentally, in most cases this increases the speed of
computation.

Certainly, one has to run the algorithm on a machine with more
number of processors to investigate its parallelization
experimentally. Because of the result of two-processor
benchmarking and the intrinsic parallel structure of the
algorithm, one can expect its good experimental scalability.

\section{Conclusions}

As discussed above, experimental analysis of the described
involutive algorithms implemented for Janet division shows that they
form an efficient computational alternative to Buchberger's algorithm
for construction of \Gr bases.

Evidently, there are infinitely many different and algorithmically
``good'', i.e. Noetherian, continuous and constructive,
involutive divisions~\cite{G02} satisfying properties in Definition 3.1.
However it is still an open question, if there is a better division in theory
and/or heuristically better in practice than Janet division. In theory,
one can say that division ${\cL}_1$ is better than ${\cL}_2$ if the inclusion
\begin{equation}
(\forall U\in \M)\ (\forall u\in U)\ [\ M_{{\cL}_1}\supseteq M_{{\cL}_2}\ ] \label{comp_div}
\end{equation}
holds for conventionally autoreduced monomial sets $U$, and
there are sets $U$ for which the inclusion is strict. In this case,
by the same arguments as those used in~\cite{GB98b} for the correctness proof, one can
show that, for an identical selection strategy, algorithm {\bf InvolutiveBasis II}
with division ${\cL}_1$ will process never more, but sometimes less number of
non-multiplicative prolongations than with division ${\cL}_2$. In addition,
the corresponding output ${\cL}-$bases satisfy $G_1\subseteq G_2$ with the strict inclusion
for some ideals.

In practice, running time depends not only on the number of
prolongation treated and on the related size of the intermediate and final
basis. Apart from other implementation aspects, running time depends also on
costs of such frequently used and depending on the involutive division operations
as re-computation of multiplicative and non-multiplicative variables under
change of the intermediate polynomial set and search for involutive divisors.
However, it is fairly advisable to prefer (for implementation too) divisions
which generate less number of non-multiplicative variables.

To illustrate the last statement, confront Pommaret and Janet divisions for
zero-dimen\-sional ideals when Pommaret bases always exist. In such a case
at first glance Pommaret division (Example 3.2) is favoured over Janet one (Example 3.1).
Indeed, it does not require re-computing separation of variables for elements in the
intermediate basis if other elements are added or removed. Besides, Pommaret division
admits trees for searching involutive divisors that are very similar to Janet trees,
and the search in a Pommaret tree is as fast as in the Janet tree with the same leaves.
Nevertheless, as implemented in the form of algorithm {\bf InvolutiveBasis II},
Pommaret division compares unfavourably with Janet division. We experimented with both
divisions while creating the C code presented in~\cite{GBY01}, and for all the
benchmarks used Pommaret division lost in speed to Janet division.

The reason is that Janet division surpasses Pommaret division in accordance to
(\ref{comp_div}) (see \cite{GB98a,G00}). Table 1 gives an example of the strict inclusion.
To provide existence of the related Pommaret basis, $U$ can be considered as the leading
monomial set for zero-dimensional ideal, for instance,
$\Id(x_1^2x_3-1,x_1x_2-1,x_1x_3^2-1)$. Since, as shown in~\cite{G00}, whenever Pommaret
basis exists it coincides with (minimal) Janet basis, algorithm {\bf InvolutiveBasis II}
for Janet division produces the same output as for Pommaret division but with processing
less number of prolongations.

For the time being we do not know involutive divisions better in theory
and/or heuristically in practice than Janet division. In the approach of
Apel~\cite{Apel} involutive division is defined locally in terms of (admissible for)
a monomial set as satisfying properties 1 and 2 in Definition 3.1 and not necessarily
property 3. In this approach, given a monomial set $U$, one can always construct a division
which surpasses Janet division for $U$ in the sense of relation (\ref{comp_div}) and
often the strict inclusion holds~\cite{Apel,Apel03}. Based on this (local) concept
of involutive division, Apel designed an algorithm for completion to involution with
dynamical construction of the best involutive division for every intermediate basis.
This, theoretically very attractive, algorithmic procedure, is unlikely practical
because of high computational costs of constructing the best intermediate divisions.
A step forward to practicality of this approach was made by Hemmecke in the form of
the SlicedDivision algorithm implemented in Aldor~\cite{Hem03,HemThesis}. Though this
algorithm is far from being competitive in efficiency with algorithm
{\bf InvolutiveBasis II}, the both approaches to theory of involutive division
are complementary~\cite{Apel03} and their further development mutually helpful.

One should also say that algorithm {\bf InvolutiveBasis II} is apparently much less
efficient than the fastest modern-day algorithms designed by Faug\`ere: $F_4$~\cite{F4}
recently built-in Magma (version 2-11) and $F_5$~\cite{F5}. These two algorithms are
based on entirely different completion strategy and exploit the linear
algebra methods for polynomial reduction rather than elementary reduction chains
used by Buchberger's algorithm and the involutive algorithm. We have good reason to
think that the algorithmic ideas of Faug\`ere can be
incorporated into the involutive methods too.

\section{Acknowledgements}
I am grateful to Denis Yanovich for his assistance in computer experiments. The research
presented in this paper was partially supported by the grants 04-01-00784 from the Russian
Foundation for Basic Research and 2339.2003.2 from the Russian
Ministry of Science and Education.


\begin{thebibliography}{99}

\bibitem{BB65} B.Buchberger, An Algorithm for Finding a Basis for the Residue Class Ring
of a Zero-Dimensional Polynomial Ideal. PhD Thesis, University of Innsbruck, Institute
for Mathematics, 1965 (in German).

\bibitem{BW} T.Becker and V.Weispfenning. \Gr Bases. A Computational Approach to
 Commutative Algebra. Graduate Texts in Mathematics 141, Springer, New York, 1993.

\bibitem{Mishra} B.Mishra, Algorithmic Algebra, Springer, New York, 1993.

\bibitem{AL} W.W.Adams and P.Loustaunau, An Introduction to Grobner Bases. Graduate Studies
in Mathematics, Vol.3, American Mathematical Society, 1994.

\bibitem{CLO1} D.Cox, J.Little and D.O'Shea, Ideals, Varieties and Algorithms. An
 Introduction to Computational Algebraic Geometry and Commutative Algebra, 2nd Edition,
 Springer, New-York, 1996.

\bibitem{CLO2} D.Cox, J.Little and D.O'Shea, Using Algebraic Geometry. Graduate Texts in
Mathematics 185, Springer, New-York, 1998.

\bibitem{KR} M.Kreutzer and L.Robbiano, Computational Commutative Algebra 1, Springer,
Berlin, 2000.

\bibitem{GP} G.-M.Greuel and G.Pfister, A Singular Introduction to Commutative Algebra,
Springer, Berlin, 2002.

\bibitem{BW98} \Gr Bases and Applications, B.Buchberger and F.Winkler, eds.,
Cambridge Universite Press, 1998.

\bibitem{CAHandbook} Computer Algebra Handbook. J.Grabmeier, E.Kaltofen and V.Weispfenning, eds.,
Springer, berlin, 2003.

\bibitem{GB98a} V.P.Gerdt and Yu.A.Blinkov, Involutive Bases of
 Polynomial Ideals, Mathematics and Computers in Simulation, 45 (1998), 519--542.

\bibitem{ZB} A.Yu.Zharkov and Yu.A.Blinkov, Involutive Approach
 to Investigating Polynomial Systems, Mathematics and Computers in Simulation 42 (1996),
 323--332.

\bibitem{Buch85} B.Buchberger, Gr\"obner Bases: an Algorithmic
 Method in Polynomial Ideal Theory, Recent Trends in
 Multidimensional System Theory, N.K. Bose (ed.), Reidel, Dordrecht
 (1985) pp. 184--232.

\bibitem{GB98b} V.P.Gerdt and Yu.A.Blinkov, Minimal Involutive Bases,
Mathematics and Computers in Simulation, 45 (1998), 543--560.

\bibitem{G99} V.P.Gerdt, Completion of Linear Differential Systems to
 Involution, Computer Algebra in Scientific Computing / CASC'99,
 V.G. Ganzha, E.W. Mayr and E.V. Vorozhtsov (eds.), Springer, Berlin,
 1999, pp.115--137.

\bibitem{GBY01} V.P.Gerdt, Yu.A.Blinkov and D.A.Yanovich,
 Construction of Janet bases. I. Monomial bases, Computer Algebra in
 Scientific Computing / CASC'01, V.G.Ganzha, E.W.Mayr and E.V.Vorozhtsov (eds.),
 Springer, Berlin, 2001, pp.233--247; II. Polynomial bases, ibid., pp.249--263.

\bibitem{G02} V.P.Gerdt, Involutive Division Technique: Some Generalizations and
 Optimizations, Journal of Mathematical Sciences 108(6) (2002), 1034-1051.

\bibitem{CG01} Y.F.Chen and X.S.Gao, Involutive Directions and New Involutive Divisions,
Computers and Mathematics with Applications, 41(7-8) (2001), 945--956.

\bibitem{CHS} J.Calmet, M.Hausdorf and W.M.Seiler, A Constructive Introduction to Involution
Proceedings of the International Symposium on Applications of Computer Algebra - ISACA 2000,
R.Akerkar, ed., Allied Publishers, New Delhi, 2001, pp.33--50.

\bibitem{S-1} M.Hausdorf, W.M.Seiler and R.Steinwandt, Involutive Bases in the Weyl Algebra,
Journal of Symbolic Computation 34 (2002), 181--198.

\bibitem{S-2} W.M.Seiler, Involution - The formal theory of differential equations and
its applications in computer algebra and numerical analysis, Habilitation thesis,
Dept. of Mathematics, University of Mannheim, 2002.

\bibitem{Semenov02} A.S.Semenov, Static Properties of Involutive Divisions, Proceedings of the
Workshop on Under-and Overdetermined Systems of Algebraic or Differential Equations, J.Calmet,
M.Hausdorf and W.M.Seiler, eds., Institute f\"ur Algorithmen and Cognitive Systeme,
University Karlsruhe, 2002, pp.151-155.

\bibitem{Carra-1} V.Marotta and G.Carra-Ferro, Involutive Division and Involutive Autoreduction,
Proceedings of the 8th Rhine Workshop on Computer Algebra, H.Kredek and W.K.Seiler, eds., University
of Mannheim, 2002, pp.115--124.

\bibitem{Carra-2} G.Carra-Ferro, M.D'Anna and V.Marotta, Characterization of Involutive Divisions and Its Applications,
Proceedings of the Workshop on Under-and Overdetermined Systems of Algebraic or Differential Equations, J.Calmet,
M.Hausdorf and W.M.Seiler, eds., Institute f\"ur Algorithmen and Cognitive Systeme,
University Karlsruhe, 2002, pp.19--36.

\bibitem{CG03} Y.Chen and X.-S. Gao, Involutive Bases of Algebraic Partial Differential
Equation Systems, Science in China (A), 33(2) (2003), 97--113.

\bibitem{Evans} G.A.Evans, Noncommutative Involutive Bases, Proceedings of ACA 04, the 10th
International Conference on Applications of Computer Algebra (July 21-23, 2004, Beaumont,
Texas, U.S.A.), to appear.

\bibitem{Berth} M.Berth and V.Gerdt, Computation of Involutive Bases with Mathematica, Proceedings
of the Third International Workshop on ``Mathematica'' System in Teaching and Research,
University of Podlasie, Seldce, Poland, 2001, pp.29--34.

\bibitem{Daniel} Yu.A.Blinkov, V.P.Gerdt, C.F.Cid, W.Plesken and D.Robertz.
 The Maple Package "Janet": I.Polynomial Systems and II.Linear
 Partial Differential Equations, Computer Algebra in Scientific
 Computing / CASC 2003, V.G.Ganzha, E.W.Mayr, and E.V.Vorozhtsov, eds., Institute of
 Informatics, Technical University of Munich, Garching, 2003, pp.31--54.; II.Linear
 Partial Differential Equations, ibid., pp.41--54.

\bibitem{SH} M.Hausdorf and W.M.Seiler, Involutive Bases in MuPAD I:
 Involutive Divisions, mathPAD 11 (2002), 51-56; II Polynomial Algebras of Solvable Type,
 mathPAD (to apper).

\bibitem{Apel} J.Apel, Theory of Involutive Divisions and an Application to Hilbert Function
Computations, Journal of Symbolic Computation, 25 (1998), 683--704.

\bibitem{Hem02} R.Hemmecke, An Efficient method for Finding a Simplifier in an Involutive Basis
Computation, Proceedings of the Workshop on Under-and Overdetermined Systems of Algebraic or
Differential Equations, J.Calmet,M.Hausdorf and W.M.Seiler, eds., Institute f\"ur Algorithmen
and Cognitive Systeme, University Karlsruhe, 2002, pp.87-96.

\bibitem{AH02} J.Apel and R.Hemmecke, Detecting unnecessary reductions in an involutive basis
computation, RISC Linz Report Series 02-22, 2002.

\bibitem{Apel03} J.Apel, Passive Complete Orthonomic Systems and Involutive Bases, Symbolic and Numeric
 Scientific Computation, F.Winkler and U.Langer, eds., Lecture Notes in Computer Science
 2630, 2003, pp.88--107.

\bibitem{Hem03} R.Hemmecke, Dynamical Aspects of Involutive Bases Computations, Symbolic and Numeric
 Scientific Computation, F. Winkler and U.Langer, eds., Lecture Notes in Computer Science
 2630, 2003, pp.168--182.

\bibitem{HemThesis} R.Hemmecke, Involutive Bases for Polynomial Ideals, PhD Thesis, RISC Linz, 2003.

\bibitem{BM} D.Bini and B.Mourrain, Polynomial Test Suite (1996).

 {\tt http://www-sop.inria.fr/saga/POL}

\bibitem{Verschelde} J.Verschelde, The Database with Test Examples.

 {\tt http://www.math.uic.edu/\~\,jan/demo.html}


\bibitem{Janet} M.Janet, Le\c cons sur les Syst\`{e}mes
 d'Equations aux D\'eriv\'ees Partielles, Cahiers Scientifiques, IV,
 Gauthier-Villars, Paris, 1929.

\bibitem{Thomas} J.Thomas, Differential Systems, American
 Mathematical Society, New York, 1937.

\bibitem{Pommaret} J.F.Pommaret, Systems of Partial
 Differential Equations and Lie Pseu\-do\-groups, Gordon \& Breach, New York, 1978.

\bibitem{G02-1} V.P.Gerdt, On an Algorithmic Optimization in Computation of
 Involutive Bases, Programming and Computer Software 28, 2 (2002), 62--65.

\bibitem{Arnold} E.A.Arnold, Modular algorithms for computing \Gr bases,
Journal of Symbolic Computation 35 (2003), 403--420.

\bibitem{Traverso} A.Giovinni, T.Mora, G.Niesi, L.Robbiano and C.Traverso,
One sugar cube, please, or selection strategies in the Buchberger algorithm,
Proceedings of ISSAC'91, ACM Press, New York, 1991, pp.49--54.

\bibitem{Blinkov} Yu.A.Blinkov, Method of Separative Monomials for Involutive Divisions,
Programming and Computer Software 27, 3 (2001), 139--141.

\bibitem{GY04} V.P.Gerdt and D.A.Yanovich, Parallel Computation of Involutive and \Gr Bases,
Computer Algebra in Scientific Computing / CASC 2003, V.G.Ganzha, E.W.Mayr, and E.V.Vorozhtsov,
eds., Institute of Informatics, Technical University of Munich, Garching, 2004, pp.185--194.

\bibitem{MTG} B.Amrheim, O.Gloor and W.K\"{u}chlin, A case study of
multi-threaded Gr\"{o}bner basis completion. Proceedings
of ISSAC'96, ACM Press, New York, 1996, pp. 95--102.

\bibitem{G00} V.P.Gerdt, On the Relation between Pommaret and Janet Bases,
 Computer Algebra in Scientific Computing / CASC'00,
 V.G. Ganzha, E.W. Mayr and E.V. Vorozhtsov, eds., Springer, Berlin,
 2000, pp.167--181.

\bibitem{F4} J.C.Faug\`{e}re, A new efficient algorithm for computing \Gr bases ($F_4$),
Journal of Pure and Applied Algebra 139, 1-3 (1999), 61--68.

\bibitem{F5} J.C.Faug\`{e}re, A new efficient algorithm for computing \Gr bases without
reduction to zero ($F_5$), Proceedings of Issac 2002, ACM Press, New York, 2002, pp. 75--83.

\end{thebibliography}
\end{document}